\documentclass[12pt,reqno]{amsart}

\usepackage{hyperref}
\usepackage{amssymb}
\usepackage{graphicx}
\usepackage{float}
\usepackage{tikz}
\usetikzlibrary{arrows}
\usepackage{pst-node}
\usepackage{tikz-cd}

\numberwithin{equation}{section}

\newtheorem{theorem}{Theorem}[section]
\newtheorem{definition}[theorem]{Definition}

\newtheorem{lemma}[theorem]{Lemma}
\newtheorem{proposition}[theorem]{Proposition}

\newtheorem{remark}[theorem]{Remark}

\newtheorem{corollary}[theorem]{Corollary}

\newtheorem*{notation*}{Notation}

\title[\textsc{Cohen--Macaulayness and Buchsbaumness}]
{Cohen--Macaulayness and Buchsbaumness of Cartesian
Products of Paths and Cycles}

\author{Bidwan Chakraborty}

\address{Department of Mathematics, National Institute of Technology Puducherry, Karaikal 609609, Puducherry, India}

\email{bidwan.c@nitpy.ac.in}

\thanks{ORCID: \href{https://orcid.org/0000-0002-7471-3359}
{0000-0002-7471-3359}}



\begin{document}

\maketitle

% Explicitly display affiliation and contact details
% on the first page.

\begin{center}
\small
Department of Mathematics\\
National Institute of Technology Puducherry\\

Email: \href{mailto:bidwan.c@nitpy.ac.in}
{bidwan.c@nitpy.ac.in}\\[4pt]
% Uncomment and replace with your verified ORCID:
% ORCID: \href{https://orcid.org/YOUR-ORCID-ID}
% {YOUR-ORCID-ID}
\end{center}

%------------------------------------------------
% Continue with your existing abstract here.
% Do not delete the mathematical content.
%------------------------------------------------

% Place your abstract and the rest of your manuscript here.

\begin{abstract}
We study the independence complexes of Cartesian products of paths and cycles. Let $G_{m,n}=P_m\square C_n$ be the Cartesian product and let $\Delta_{m,n}=\operatorname{Ind}(G_{m,n})$ be its independence complex. We determine the well-covered graphs in this family and use reduced homology to study their Cohen–Macaulay and Buchsbaum properties. Using these results, we completely classify the Cohen--Macaulay members. We then study the Buchsbaum property of $\Delta_{m,n}$ through vertex links. For the family $P_m\square C_3$, we introduce an auxiliary family of graphs and use an inductive link argument together with Reisner's criterion to determine precisely when $\operatorname{Ind}(P_m\square C_3)$ is Buchsbaum. Combining this with the remaining well-covered cases yields the complete Buchsbaum classification.
\end{abstract}

\section{Introduction}

Let $G=(V_G,E_G)$ be a finite simple graph with vertex set
$V_G=\{x_1,\ldots,x_n\}$, and let $R=k[x_1,\ldots,x_n]$, where $k$ is a
field. The edge ideal of $G$ is the squarefree quadratic monomial ideal
$
I(G)=(x_ix_j:\{x_i,x_j\}\in E_G)
$
in $R$. Since $I(G)$ is the Stanley-Reisner ideal of the independence
complex $\operatorname{Ind}(G)$, the quotient ring $R/I(G)$ is
Cohen--Macaulay if and only if $\operatorname{Ind}(G)$ is
Cohen--Macaulay. Thus, the study of Cohen--Macaulay graphs provides a
natural connection between graph theory, combinatorial topology, and
commutative algebra. In particular, Reisner's criterion characterizes
the Cohen--Macaulay property of a simplicial complex in terms of the
reduced homology of its links.
For background on Cohen--Macaulay graphs and related
classes, see, for example,
\cite{CookNagel2012,DaoDoolittleLyle2020,
earl2016independence,EstradaVillarreal1997,
FranciscoHa2009,HerzogHibiZheng2006,
Stanley1977,Villarreal1990}.

A closely related notion is that of a Buchsbaum simplicial complex.
A pure simplicial complex $\Delta$ is Buchsbaum over $k$ if and only if
$\operatorname{lk}_{\Delta}(v)$ is Cohen--Macaulay over $k$ for every
vertex $v\in V(\Delta)$. Every Cohen--Macaulay simplicial complex is
Buchsbaum. Consequently, for independence complexes, the study of
vertex links provides a useful method for determining the Buchsbaum
property. Cook and Nagel \cite{CookNagel2012} investigated several
characterizations of graphs whose independence complexes are
Cohen--Macaulay or Buchsbaum, including connections with flag complexes
and vertex-decomposability. Buchsbaum properties of independence
complexes of well-covered graphs have also been studied in
\cite{earl2016independence}.

Independence complexes of graph products provide a natural setting for
investigating these properties. Let $G$ and $H$ be graphs. The \emph{Cartesian product}
$G\square H$ is the graph with vertex set
$V(G\square H)=V(G)\times V(H)$, where two vertices $(g,h)$ and
$(g',h')$ are adjacent if and only if either $g=g'$ and
$\{h,h'\}\in E(H)$, or $h=h'$ and $\{g,g'\}\in E(G)$. Hamming graphs,
which are Cartesian products of complete graphs, were investigated in
\cite{rao2026cohen}. The Cartesian product of a path with a cycle has a
different local structure and provides a natural family for studying how
graph geometry influences the Cohen--Macaulay and Buchsbaum properties
of independence complexes.

For an independent set $\sigma$ of a graph $G$, the link in its
independence complex satisfies
$
\operatorname{lk}_{\operatorname{Ind}(G)}(\sigma)
=
\operatorname{Ind}(G-N_G[\sigma]).
$
Thus, the study of links in an independence complex can be reduced to
the study of independence complexes of graphs obtained by deleting
closed neighborhoods. We will use this observation repeatedly
throughout the paper. We also use the fact that a graph $G$ is
well-covered if and only if $\operatorname{Ind}(G)$ is pure.

For $m\geq1$, let $P_m$ denote the path graph on vertices
$1,\ldots,m$ with edges $\{i,i+1\}$ for $1\leq i\leq m-1$, and for
$n\geq3$, let $C_n$ denote the cycle graph on $\mathbb{Z}_n$ with edges
$\{i,i+1\}$, where the indices are taken modulo $n$. We consider the
Cartesian product
$
G_{m,n}=P_m\square C_n,
$
whose vertex set is
$
V(G_{m,n})=\{(i,j):1\leq i\leq m,\ j\in\mathbb{Z}_n\},
$
and write
$
\Delta_{m,n}=\operatorname{Ind}(G_{m,n}).
$
Our first step is to determine precisely when $G_{m,n}$ is well-covered.
 Using a result of Hartnell, Rall, and Wash on well-covered
Cartesian products, together with a direct analysis of the remaining
cases, we obtain the following classification.

\vskip 2mm \noindent
{\bf Theorem \ref{thm:main-well-covered}.}
[Classification of well-covered products]

Let $m\geq 1$ and $n\geq 3$. Then $P_m\square C_n$ is well-covered if
and only if
$
n=3
\quad\text{or}\quad
(m,n)\in\{(1,4),(1,5),(1,7),(2,5)\}.
$
Equivalently, $P_m\square C_n$ is well-covered precisely in the following
cases:
\vspace{-.1cm}
\begin{enumerate}
    \item $n=3$ and $m\geq 1$;
    \item $m=1$ and $n\in\{4,5,7\}$;
    \item $m=2$ and $n=5$.
\end{enumerate}
Using this, we proceed to study the Cohen–Macaulay and Buchsbaum properties, since both properties imply purity of the independence
complex. Theorem \ref{thm:main-well-covered} reduces the study of Cohen--Macaulayness and Buchsbaumness to the surviving cases
$P_m\square C_3$, $C_4$, $C_5$, $C_7$, and $P_2\square C_5$.

The case $P_m\square C_3$ requires a separate analysis. The graph
$P_m\square C_3$ consists of $m$ consecutive columns, each inducing a
triangle, with matching edges joining corresponding vertices in
consecutive columns. This triangular-column structure makes the case
$n=3$ substantially different from the cases $n\geq4$. In particular,
it allows vertex links near the boundary to be described explicitly in
terms of independence complexes of smaller auxiliary graphs.

For this purpose, we introduce a family of graphs $T_k$, consisting of
an initial edge followed by $k-1$ triangular columns. We use these graphs
to describe the links of boundary vertices in
$\operatorname{Ind}(P_m\square C_3)$. We prove that
$\operatorname{Ind}(T_k)$ is not Cohen--Macaulay for every $k\geq3$.
Consequently, for $m\geq4$, a boundary vertex of $P_m\square C_3$ has a
link that is not Cohen--Macaulay, which gives the failure of the
Buchsbaum property. The cases $m=1,2,3$ are treated separately.
We next determine exactly which of these independence complexes are
Cohen--Macaulay.

\vskip 2mm \noindent
{\bf Theorem \ref{CompleteCM}.}
The independence complex
$\operatorname{Ind}(P_m\square C_n)$ is Cohen--Macaulay if and only if
$
(m,n)\in\{(1,3),(2,3),(1,5)\}.
$

Combining this analysis with the remaining well-covered cases gives the
complete Buchsbaum classification.

\vskip 2mm \noindent
{\bf Theorem \ref{thm:complete-buchsbaum}.}
Let $m\geq1$ and $n\geq3$. Then $G_{m,n}=P_m\square C_n$ is Buchsbaum if
and only if
$
n=3\text{ and }m\in\{1,2,3\},
$
or
$
(m,n)\in\{(1,4),(1,5),(2,5),(1,7)\}.
$

Thus, although every Cohen--Macaulay independence complex is Buchsbaum,
the two properties differ precisely for
$
(m,n)\in\{(1,4),(1,7),(2,5),(3,3)\}.
$
This distinction is governed by the local structure of the Cartesian
product, with the case $C_3=K_3$ playing a particularly important role.
The auxiliary graphs $T_k$ provide an explicit description of the
obstructions to Buchsbaumness in this case.

The paper is organized as follows. In Section 2, we recall the necessary
definitions and results concerning independence complexes, links,
Cohen--Macaulay and Buchsbaum complexes, and Reisner's criterion. In
Section 3, we classify the well-covered graphs $P_m \square C_n$ and
determine the Cohen--Macaulay cases. In Section 4, we introduce the
auxiliary graphs $T_k$, determine the Buchsbaum property for
$P_m \square C_3$, and combine this with the remaining well-covered
cases to obtain the complete Cohen--Macaulay and Buchsbaum classification.

\section{Preliminaries}
In this section, we review the relevant definitions and results.

We recall the basic notions from simplicial complex theory that will be used throughout the paper. A simplicial complex $\Delta$ on a finite vertex set $V(\Delta)$ is a collection of subsets of $V(\Delta)$ containing the empty set and closed under taking subsets. A face $F$ of $\Delta$ is called a facet if it is maximal with respect to inclusion, and its dimension is defined by $\dim F=|F|-1$. The dimension of $\Delta$ is the maximum dimension of its faces. The complex $\Delta$ is called pure if all its facets have the same dimension. For a face $\sigma\in\Delta$, the link of $\sigma$ in $\Delta$ is the complex $\operatorname{lk}_{\Delta}(\sigma)=\{\tau\in\Delta:\tau\cap\sigma=\emptyset,\ \tau\cup\sigma\in\Delta\}$, while the deletion of a vertex $v\in V(\Delta)$ is $\operatorname{del}_{\Delta}(v)=\{\sigma\in\Delta:v\notin\sigma\}$. If $\Delta$ and $\Gamma$ have disjoint vertex sets, their join is the complex $\Delta*\Gamma=\{\sigma\cup\tau:\sigma\in\Delta,\ \tau\in\Gamma\}$. The reduced homology groups of $\Delta$ over a field $k$ are denoted by $\widetilde{H}_i(\Delta;k)$. A simplicial complex $\Delta$ is Cohen--Macaulay over $k$ if it is pure and $\widetilde{H}_i(\operatorname{lk}_{\Delta}(\sigma);k)=0$ for every $\sigma\in\Delta$ and every $i<\dim(\operatorname{lk}_{\Delta}(\sigma))$. A pure simplicial complex is Buchsbaum over $k$ if the link of every vertex is Cohen--Macaulay over $k$. We also use the notions of shellability and vertex-decomposability in their standard senses: $\Delta$ is shellable if its facets admit an ordering $F_1,\ldots,F_t$ such that, for every $j>1$, the complex generated by $F_j\cap F_1,\ldots,F_j\cap F_{j-1}$ is pure of dimension $\dim F_j-1$; and $\Delta$ is vertex-decomposable if it is a simplex or there is a vertex $v$ such that $\operatorname{del}_{\Delta}(v)$ and $\operatorname{lk}_{\Delta}(v)$ are vertex-decomposable and every facet of $\operatorname{del}_{\Delta}(v)$ is a facet of $\Delta$. For a simplicial complex $\Delta$ on $[n]$, its Stanley--Reisner ideal is $I_{\Delta}=(x_F:F\subseteq[n],\ F\notin\Delta)\subseteq k[x_1,\ldots,x_n]$, where $x_F=\prod_{i\in F}x_i$, and its Stanley--Reisner ring is $k[\Delta]=k[x_1,\ldots,x_n]/I_{\Delta}$. These notions provide the simplicial-complex framework for the graph-theoretic properties considered below.

We now recall the graph-theoretic terminology used in the paper. For a graph $G$, the open neighborhood of a vertex $v$ is $N_G(v)=\{u\in V(G):\{u,v\}\in E(G)\}$, and its closed neighborhood is $N_G[v]=N_G(v)\cup\{v\}$. For an independent set $\sigma$ of $G$, we write $N_G[\sigma]=\bigcup_{v\in\sigma}N_G[v]$. The girth of $G$, denoted by $g(G)$, is the length of a shortest cycle, with $g(G)=\infty$ when $G$ is acyclic. An independent set of $G$ is a subset of $V(G)$ containing no two adjacent vertices, and the independence number $\alpha(G)$ is the maximum cardinality of an independent set. An independent set is maximal if it is not properly contained in another independent set, and the independent domination number $i(G)$ is the minimum cardinality of a maximal independent set. Thus, $G$ is well-covered if all its maximal independent sets have the same cardinality, equivalently, if $i(G)=\alpha(G)$. The independence complex of $G$, denoted by $\operatorname{Ind}(G)$, is the simplicial complex whose faces are precisely the independent sets of $G$. In particular, $G$ is well-covered if and only if $\operatorname{Ind}(G)$ is pure.

\section{Well-coveredness of Cartesian products of paths and cycles}

Let $G_{m,n}=P_m\square C_n$, where $m\geq 1$ and $n\geq 3$. We determine precisely when $G_{m,n}$ is well-covered.

We shall use the following result on Cartesian products.

\begin{theorem}[\cite{HartnellRallWash2018}]
\label{thm:HRW}
Let $G$ and $H$ be nontrivial connected graphs, both of girth at least
$4$. If $G\square H$ is well-covered, then one of $G$ and $H$ is
isomorphic to $K_2$. Moreover, if both $G$ and $H$ have girth at least
$5$, then $G\square H$ is well-covered if and only if
$G\square H\cong K_2\square K_2$ or
$G\square H\cong K_2\square C_5$.
\end{theorem}

\begin{remark}
The first assertion of Theorem~\ref{thm:HRW} will be used to eliminate
all products $P_m\square C_n$ with $m\geq 3$ and $n\geq 4$. The second
assertion will be used for the case $m=2$.
\end{remark}

We first deal with the case $n=3$.

\begin{lemma}
\label{lem:C3-well-covered}
For every $m\geq 1$, the graph $P_m\square C_3$ is well-covered.
\end{lemma}

\begin{proof}
Let us write the vertices of $P_m\square C_3$ as
$V(P_m\square C_3)=\{v_{i,j}:1\leq i\leq m,\ 0\leq j\leq 2\}$.
For each $i$, let $T_i=\{v_{i,0},v_{i,1},v_{i,2}\}$. The subgraph
induced by $T_i$ is a triangle. Consequently, every independent set
contains at most one vertex from each $T_i$. We claim that every maximal independent set contains exactly one vertex
from every $T_i$. Suppose, to the contrary, that $M$ is a maximal
independent set and that $T_i\cap M=\varnothing$ for some $i$. If $1<i<m$, the vertices of $T_i$ which can be prevented from being
added to $M$ are adjacent only to vertices belonging to the two
neighboring triangles $T_{i-1}$ and $T_{i+1}$. Since $M$ contains at
most one vertex from each of these triangles, at most two vertices of
$T_i$ are blocked. Hence, at least one vertex of $T_i$ can be added to
$M$, contradicting the maximality of $M$.

If $i=1$ or $i=m$, there is only one neighboring triangle, so at most
one vertex of $T_i$ is blocked. Again, this contradicts the maximality
of $M$. Therefore, every maximal independent set intersects every $T_i$. Since
an independent set contains at most one vertex from each $T_i$, every
maximal independent set has exactly $m$ vertices. Hence, $i(P_m\square C_3)=\alpha(P_m\square C_3)=m$, and therefore
$P_m\square C_3$ is well-covered.
\end{proof}

\begin{remark}\label{well-cycle}

The case $m=1$ reduces to cycles. It is well known that
$i(C_n)=\left\lceil n/3\right\rceil$ and
$\alpha(C_n)=\left\lfloor n/2\right\rfloor$.
Hence, $C_n$ is well-covered if and only if
$\left\lfloor n/2\right\rfloor=\left\lceil n/3\right\rceil$.
A direct verification gives $n=3,4,5,7$. For $n\geq 8$, we have $\left\lfloor n/2\right\rfloor>
\left\lceil n/3\right\rceil$, and consequently $C_n$ is not well-covered. Therefore, $C_n$ is well-covered if and only if $n\in\{3,4,5,7\}$.
\end{remark}

We now exclude all remaining cases with $m\geq 3$.

\begin{lemma}
\label{lem:mgeq3}
Let $m\geq 3$ and $n\geq 4$. Then $P_m\square C_n$ is not well-covered.
\end{lemma}

\begin{proof}
The graphs $P_m$ and $C_n$ are nontrivial connected graphs. Moreover,
$\operatorname{girth}(P_m)=\infty$ and
$\operatorname{girth}(C_n)=n\geq 4$. Thus both factors have girth at
least $4$. If $P_m\square C_n$ were well-covered, Theorem~\ref{thm:HRW} would
imply that one of the factors is isomorphic to $K_2$. However,
$P_m\not\cong K_2$ for $m\geq 3$, and
$C_n\not\cong K_2$ for $n\geq 4$. This is a contradiction.
Therefore, $P_m\square C_n$ is not well-covered.
\end{proof}

It remains to consider the case $m=2$.

\begin{lemma}
\label{lem:m2}
For $n\geq 4$, the graph $P_2\square C_n$ is well-covered if and only if
$n=5$.
\end{lemma}

\begin{proof}
Since $P_2\cong K_2$, we consider two cases. First, let $n=4$. Label the vertices of $P_2\square C_4$ by
$(i,j)$, where $i\in\{1,2\}$ and $j\in\mathbb{Z}_4$. Then
$M_1=\{(1,0),(2,2)\}$ is a maximal independent set of cardinality $2$,
whereas $M_2=\{(1,0),(1,2),(2,1),(2,3)\}$ is a maximal independent set of cardinality $4$. Hence, $P_2\square C_4$ is not well-covered.

Now suppose that $n\geq 5$. Both $K_2$ and $C_n$ are nontrivial
connected graphs of girth at least $5$. By the second assertion of
Theorem~\ref{thm:HRW}, a well-covered Cartesian product of two such
graphs can only be $K_2\square K_2$ or $K_2\square C_5$.
Consequently, $K_2\square C_n$ is well-covered only when $n=5$.
Conversely, $K_2\square C_5$ is well-covered by the same theorem.
Therefore
$P_2\square C_n$ is well-covered if and only if $n=5$.
\end{proof}

We can now state the complete classification.

\begin{theorem}[Classification of well-covered products]
\label{thm:main-well-covered}
Let $m\geq 1$ and $n\geq 3$. Then $P_m\square C_n$ is well-covered if
and only if
$n=3$ or
$(m,n)\in\{(1,4),(1,5),(1,7),(2,5)\}$.
Equivalently,
$P_m\square C_n$ is well-covered precisely in the following cases:
\begin{enumerate}
    \item $n=3$ and $m\geq 1$;
    \item $m=1$ and $n\in\{4,5,7\}$;
    \item $m=2$ and $n=5$.
\end{enumerate}
\end{theorem}

\begin{proof}
If $n=3$, the assertion follows from Lemma \ref{lem:C3-well-covered}
. If $m=1$, then
$P_1\square C_n\cong C_n$, and the assertion follows from
Remark {\ref{well-cycle}}.

It remains to consider $m\geq 2$ and $n\geq 4$. If $m\geq 3$, the result follows from Lemma~\ref{lem:mgeq3}. If $m=2$, the result follows from Lemma~\ref{lem:m2}. Hence, the listed cases are precisely all well-covered members of the family.
\end{proof}

\begin{corollary}
\label{cor:CM-reduction}
Let $\Delta_{m,n}=\operatorname{Ind}(P_m\square C_n)$ be the independence
complex of $P_m\square C_n$. 
If $n\neq3$ and
$(m,n)\notin\{(1,4),(1,5),(1,7),(2,5)\}$, then $\Delta_{m,n}$ is not Cohen--Macaulay
\end{corollary}

\begin{proof}
If $\Delta_{m,n}$ is Cohen--Macaulay, then $\Delta_{m,n}$ is pure.
For an independence complex, purity is equivalent to the corresponding
graph being well-covered. The result therefore follows immediately
from Theorem~\ref{thm:main-well-covered}.
\end{proof}

\begin{remark}
Theorem~\ref{thm:main-well-covered} reduces the subsequent study of
 Cohen--Macaulayness and Buchsbaumness to the surviving cases
$P_m\square C_3$, $C_4$, $C_5$, $C_7$, and $P_2\square C_5$.
In particular, every Cohen--Macaulay member of the family must belong
to one of these cases.
\end{remark}

Next we investigate Cohen--Macaulayness of $\operatorname{Ind}(P_m\square C_3).$ Let $G_m=P_m\square C_3$ and let $\Delta_m=\operatorname{Ind}(G_m)$. We determine precisely when
$\Delta_m$ is Cohen--Macaulay.

We use Reisner's criterion in the following form.

\begin{theorem}[Reisner's criterion]
\label{thm:Reisner}
Let $\Delta$ be a pure simplicial complex of dimension $d$ over a
field $k$. Then $\Delta$ is Cohen--Macaulay over $k$ if and only if,
for every face $\sigma\in\Delta$,
$\widetilde{H}_i(\operatorname{lk}_{\Delta}(\sigma);k)=0$
for all $i<\dim\operatorname{lk}_{\Delta}(\sigma)$.
\end{theorem}

For an independence complex, links have a particularly useful
description.

\begin{lemma}
\label{lem:link-independence}
Let $G$ be a graph and let $\Delta=\operatorname{Ind}(G)$. For every
independent set $\sigma\in\Delta$, we have
$\operatorname{lk}_{\Delta}(\sigma)
=
\operatorname{Ind}(G-N_G[\sigma])$,
where $N_G[\sigma]$ denotes the closed neighborhood of $\sigma$.
\end{lemma}

\begin{proof}
A face $\tau$ belongs to $\operatorname{lk}_{\Delta}(\sigma)$ if and
only if $\tau\cap\sigma=\varnothing$ and $\tau\cup\sigma$ is an
independent set of $G$. This is equivalent to saying that no vertex of
$\tau$ belongs to $N_G[\sigma]$. Hence $\tau$ is an independent set of
the induced graph $G-N_G[\sigma]$, which proves the assertion.
\end{proof}

We first consider the case $m\geq4$.

\begin{lemma}
\label{lem:general-link}
For every $m\geq4$, there exists a face $\sigma_m\in\Delta_m$ such that
$\operatorname{lk}_{\Delta_m}(\sigma_m)
\cong\operatorname{Ind}(C_4)$.
\end{lemma}

\begin{proof}
Label the vertices of $G_m=P_m\square C_3$ by $(i,j)$, where
$0\leq i\leq m-1$ and $j\in\mathbb{Z}_3$. For each fixed $i$, the
vertices $(i,0),(i,1),(i,2)$ induce a triangle.

For $m\geq4$, define
$\sigma_m=
\{(i,i\bmod 2):0\leq i\leq m-4\}
\cup\{(m-1,\varepsilon_m)\}$,
where
$\varepsilon_m=0$ if $m$ is even and
$\varepsilon_m=1$ if $m$ is odd. We first verify that $\sigma_m$ is an independent set. For consecutive
vertices in the first set, the column indices differ by one, whereas
their second coordinates have opposite parity. Hence, they are not
adjacent in the Cartesian product. The final vertex lies in column
$m-1$, whereas the last vertex of the first set lies in column $m-4$.
Thus, these two vertices are not adjacent. Hence $\sigma_m$ is
independent. We now determine the graph obtained after deleting the closed
neighborhood of $\sigma_m$. The selected vertices in consecutive columns block the corresponding vertices in the neighboring columns. Since each column induces a
triangle, the selected vertex in each column also blocks the other two vertices in that column. Consequently, all columns up to and including $m-4$ are deleted.

Similarly, the selected vertex in column $m-1$ deletes the entire last
column. Its neighbor in column $m-2$ is also deleted. Thus, the only columns that may contain vertices outside
$N_{G_m}[\sigma_m]$ are columns $m-3$ and $m-2$.

In column $m-3$, the selected vertex in column $m-4$ blocks exactly one
vertex, leaving two vertices. In column $m-2$, the selected vertex in
column $m-1$ blocks exactly one vertex, leaving two vertices.

The two surviving vertices in each column are adjacent, since each
column induces a triangle. Moreover, the Cartesian-product matching
edges join the corresponding surviving vertices in the two consecutive
columns. Hence, the remaining graph is a $4$-cycle. Therefore
$G_m-N_{G_m}[\sigma_m]\cong C_4$. By Lemma~\ref{lem:link-independence}, we obtain
$\operatorname{lk}_{\Delta_m}(\sigma_m)
\cong\operatorname{Ind}(C_4)$.
\end{proof}

We next recall the independence complex of $C_4$.

\begin{lemma}
\label{lem:indC4}
The simplicial complex $\operatorname{Ind}(C_4)$ is disconnected and
has dimension $1$. In particular,
$\widetilde{H}_0(\operatorname{Ind}(C_4);k)\neq0$
for every field $k$.
\end{lemma}

\begin{proof}
Let the vertices of $C_4$ be $v_1,v_2,v_3,v_4$ in cyclic order. The
only maximal independent sets are $\{v_1,v_3\}$ and
$\{v_2,v_4\}$. Hence $\operatorname{Ind}(C_4)$ is the disjoint union
of two edges. Therefore, it is disconnected and has dimension $1$.
Consequently,
$\widetilde{H}_0(\operatorname{Ind}(C_4);k)\neq0$.
\end{proof}

We can now establish the non-Cohen--Macaulayness for all $m\geq4$.

\begin{theorem}
\label{thm:notCM-mgeq4}
For every $m\geq4$, the independence complex
$\operatorname{Ind}(P_m\square C_3)$ is not Cohen--Macaulay over any
field $k$.
\end{theorem}

\begin{proof}
Let $\Delta_m=\operatorname{Ind}(P_m\square C_3)$. By
Lemma~\ref{lem:general-link}, there exists a face $\sigma_m\in\Delta_m$
such that
$\operatorname{lk}_{\Delta_m}(\sigma_m)
\cong\operatorname{Ind}(C_4)$. By Lemma~\ref{lem:indC4},
$\widetilde{H}_0(\operatorname{lk}_{\Delta_m}(\sigma_m);k)\neq0$.
Furthermore, $\dim\operatorname{lk}_{\Delta_m}(\sigma_m)=1$. Thus, the reduced homology in degree zero does not vanish, even though
$0<\dim\operatorname{lk}_{\Delta_m}(\sigma_m)$. Hence, the condition in
Reisner's criterion fails for the face $\sigma_m$. Therefore
$\Delta_m$ is not Cohen--Macaulay.
\end{proof}

It remains to consider the case $m=3$.

\begin{lemma}\label{lem:notCM-m3}
Let $k$ be any field and let $\Delta_3=\operatorname{Ind}(P_3\square C_3)$. Then $\Delta_3$ is not Cohen--Macaulay over $k$.
\end{lemma}

\begin{proof}
Let $G=P_3\square C_3$. Let $V(G)={(i,j):0\leq i\leq2,\ j\in\mathbb{Z}_3}$. Since each column induces a triangle, every independent set contains at most one vertex from each column. Hence $\alpha(G)\leq3$. On the other hand, ${(0,0),(1,1),(2,2)}$ is an independent set. Thus $\alpha(G)=3$ and consequently $\dim\Delta_3=2$.

We next determine the $f$-vector of $\Delta_3$. Clearly, $f_{-1}=1$ and $f_0=9$. The graph $G$ has $3\binom{3}{2}+2\cdot3=15$ edges. Hence the number of independent $2$-sets is $f_1=\binom{9}{2}-15=21$.

An independent $3$-set contains exactly one vertex from each column. Write the corresponding labels as $c_0,c_1,c_2\in\mathbb{Z}_3$. The Cartesian-product matching edges impose $c_0\neq c_1$ and $c_1\neq c_2$. There are $3$ choices for $c_0$, $2$ choices for $c_1$, and $2$ choices for $c_2$. Therefore $f_2=3\cdot2\cdot2=12$.

We claim that $\Delta_3$ is pure. Let $S$ be an independent set of $G$. If $S$ contains vertices from all three columns, then $|S|=3$. Suppose that $S$ contains vertices from exactly two columns. If these columns are consecutive, say columns $0$ and $1$, then we may choose in column $2$ any vertex whose label differs from the label of the vertex in column $1$. This vertex is not adjacent to the vertex in column $0$, since columns $0$ and $2$ are not adjacent. Hence, $S$ can be extended to an independent set of size $3$. The case of columns $1$ and $2$ is analogous.

If $S$ contains vertices from columns $0$ and $2$, say $(0,a)$ and $(2,b)$, then choose a label $c\in\mathbb{Z}_3$ satisfying $c\neq a$ and $c\neq b$. Such a label exists because $\mathbb{Z}_3$ has three elements. Then $(1,c)$ is not adjacent to either $(0,a)$ or $(2,b)$, so again $S$ can be extended to an independent set of size $3$.

Finally, if $S$ contains a vertex from only one column, then $S$ can clearly be extended to an independent set containing vertices from all three columns. Therefore, every maximal independent set has cardinality $3$, and $\Delta_3$ is pure of dimension $2$.

Put $d=\dim\Delta_3+1=3$. The $h$-vector of $\Delta_3$ is defined by $\sum_{i=0}^{3}h_i t^{3-i}=\sum_{i=0}^{3}f_{i-1}(t-1)^{3-i}$. Using the above values of the $f$-vector, we obtain
$\sum_{i=0}^{3}h_i t^{3-i}=(t-1)^3+9(t-1)^2+21(t-1)+12=t^3+6t^2+6t-1$.
Consequently, $h(\Delta_3)=(1,6,6,-1)$. As, $h_3=-1<0$ so, $\operatorname{Ind}(P_3\square C_3)$ is not Cohen--Macaulay over $k$.
\end{proof}

We can now state the complete classification.

\begin{theorem}
\label{thm:CM-PmC3}
The independence complex $\operatorname{Ind}(P_m\square C_3)$ is
Cohen--Macaulay over a field $k$ if and only if $m\leq2$.
\end{theorem}

\begin{proof}
For $m=1$, we have $P_1\square C_3\cong C_3$, and
$\operatorname{Ind}(C_3)$ is a $0$-dimensional simplicial complex.
Hence, it is Cohen--Macaulay.

For $m=2$, the graph $P_2\square C_3$ is isomorphic to
$K_2\square K_3$. Its independence complex is a $6$-cycle. Hence, it is
a connected pure $1$-dimensional simplicial complex and is therefore
Cohen--Macaulay.

For $m=3$, the assertion follows from
Lemma~\ref{lem:notCM-m3}. For every $m\geq4$, it follows from
Theorem~\ref{thm:notCM-mgeq4}. Therefore
$\operatorname{Ind}(P_m\square C_3)$ is Cohen--Macaulay if and only if
$m\leq2$.
\end{proof}

\section{Buchsbaum property}

In this section, we determine the Buchsbaum property of
$\Delta_{m,n}=\operatorname{Ind}(P_m\square C_n)$. We first study the
case $P_m\square C_3$. By
introducing the auxiliary graphs $T_k$, we describe certain vertex
links and use them to determine precisely when
$\operatorname{Ind}(P_m\square C_3)$ is Buchsbaum. We then consider the remaining well-covered cases from
Theorem \ref{thm:main-well-covered}, namely $C_4$, $C_5$, $C_7$, and
$P_2\square C_5$, and determine their Buchsbaum properties. Since a
Buchsbaum complex is pure, only the well-covered cases need to be
considered. Combining these results gives the complete Buchsbaum
classification.

We shall use the following standard facts.

\begin{lemma}
Let $A$ and $B$ be vertex-disjoint graphs with no edges between them. Then
$\operatorname{Ind}(A\sqcup B)=\operatorname{Ind}(A)*\operatorname{Ind}(B)$.
\end{lemma}

\begin{proof}
A subset $F\subseteq V(A)\sqcup V(B)$ is independent in $A\sqcup B$ if and only if $F\cap V(A)$ is independent in $A$ and $F\cap V(B)$ is independent in $B$. This is precisely the definition of a face of the simplicial join
$\operatorname{Ind}(A)*\operatorname{Ind}(B)$.
\end{proof}

\begin{lemma}
A pure one-dimensional simplicial complex is Cohen--Macaulay if and only if it is connected.
\end{lemma}

\begin{proof}
Let $\Delta$ be a pure one-dimensional simplicial complex. By Reisner's criterion, the condition corresponding to the empty face is
$\widetilde H_i(\Delta;k)=0$ for all $i<1$. Since $\Delta$ is nonempty, $\widetilde H_{-1}(\Delta;k)=0$, and $\widetilde H_0(\Delta;k)=0$ if and only if $\Delta$ is connected. For any vertex $v$, purity ensures that $\operatorname{lk}_{\Delta}(v)$ is a nonempty zero-dimensional complex. Hence, its required reduced homology vanishing for $i<0$ holds automatically. For any edge $e$, the link is $\{\emptyset\}$, which has dimension $-1$ and imposes no further condition in Reisner's criterion. Therefore, $\Delta$ is Cohen--Macaulay if and only if it is connected.
\end{proof}

Now let us concentrate in the case $n=3.$ Since $C_3=K_3$, write $G_m=G_{m,3}=P_m\square C_3$. The graph $G_m$ consists of $m$ columns, each of which induces a triangle, with matching edges joining equal labels in consecutive columns.

Now we introduce the graph that occurs naturally as a vertex-link.

\begin{definition}
For $k\geq1$, let $T_k$ be the graph with $k$ columns defined as follows. The first column consists of the edge on the vertices $(0,0)$ and $(0,1)$, while every subsequent column induces a triangle on the labels $0,1,2$. Consecutive columns are joined by matching edges between equal labels.
\end{definition}

\begin{lemma}
\label{lem:boundary-link}
Let $m\geq2$ and let $x=(0,j)$ be a vertex of $G_m$ in the first column. If $\Delta_{m,3}=\operatorname{Ind}(G_m)$, then
$\operatorname{link}_{\Delta_{m,3}}(x)\cong\operatorname{Ind}(T_{m-1})$.
\end{lemma}

\begin{proof}
By symmetry of the three labels, it is enough to consider $x=(0,0)$.

The closed neighborhood of $x$ is
$N[x]={(0,0),(0,1),(0,2),(1,0)}$.
Thus, after deleting $N[x]$, the first column disappears completely, and the second column loses the vertex $(1,0)$. The remaining vertices in the second column are $(1,1)$ and $(1,2)$, which form an edge, while all subsequent columns remain triangles with the usual matching edges. After shifting the column indices and relabelling the remaining labels, the resulting graph is precisely $T_{m-1}$. Since the link of a face in an independence complex is the independence complex of the graph obtained by deleting its closed neighborhood, the desired isomorphism follows.
\end{proof}

\begin{proposition}\label{notCM}
For every $k\geq3$, the independence complex $\operatorname{Ind}(T_k)$ is not Cohen--Macaulay.
\end{proposition}

\begin{proof}
We proceed by induction on $k$. For $k=3$, consider the vertex $x=(2,2)$ of $T_3$. Since the third column of $T_3$ is a triangle, the closed neighborhood $N_{T_3}[x]$ contains all three vertices of that column. In addition, the matching edge incident with $x$ removes the vertex $(1,2)$. Consequently,
$T_3-N_{T_3}[x]={(0,0),(0,1),(1,0),(1,1)}$.
The edges of this induced subgraph are $(0,0)-(0,1)$, $(0,0)-(1,0)$, $(0,1)-(1,1)$, and $(1,0)-(1,1)$. Thus $T_3-N_{T_3}[x]\cong C_4$.

For every graph $H$ and every vertex $x$ of $H$, we have
$\operatorname{link}{\operatorname{Ind}(H)}(x)=\operatorname{Ind}(H-N_H[x])$.
Therefore,
$\operatorname{link}{\operatorname{Ind}(T_3)}(x)\cong\operatorname{Ind}(C_4)$.
The independence complex $\operatorname{Ind}(C_4)$ consists of the two diagonals of $C_4$ together with their vertices. Hence, it is the disjoint union of two edges. In particular, it is a pure one-dimensional disconnected simplicial complex. Therefore, $\widetilde{H}_0(\operatorname{Ind}(C_4);k)\neq0$. By Reisner's criterion, $\operatorname{Ind}(C_4)$ is not Cohen--Macaulay. Hence, $\operatorname{Ind}(T_3)$ is not Cohen--Macaulay.

Now assume that $k\geq4$ and that the assertion holds for all integers $3\leq r<k$. Take $x=(0,0)$, a vertex in the initial $K_2$ of $T_k$. The closed neighborhood of $x$ consists of the two vertices of the initial $K_2$ together with the vertex $(1,0)$ in the next column. Thus, deleting $N_{T_k}[x]$ removes the initial column completely and deletes one vertex from the second column. The remaining part of the second column is therefore a $K_2$, while all subsequent columns are triangles.

After shifting the column indices by one and relabelling the two surviving vertices of the first remaining column, we obtain
$T_k-N_{T_k}[x]\cong T_{k-1}$.
Consequently,
$\operatorname{link}{\operatorname{Ind}(T_k)}(x)\cong\operatorname{Ind}(T{k-1})$.
Since $k\geq4$, we have $k-1\geq3$, and hence, by the induction hypothesis, $\operatorname{Ind}(T_{k-1})$ is not Cohen--Macaulay. Therefore, $\operatorname{link}_{\operatorname{Ind}(T_k)}(x)$ is not Cohen--Macaulay. If $\operatorname{Ind}(T_k)$ were Cohen--Macaulay, then the link of every vertex would be Cohen--Macaulay. This contradiction shows that $\operatorname{Ind}(T_k)$ is not Cohen--Macaulay. Thus, by induction, $\operatorname{Ind}(T_k)$ is not Cohen--Macaulay for every $k\geq3$.
\end{proof}

We now determine the Buchsbaum property of $G_m$.

\begin{theorem}
\label{thm:C3-buchsbaum}
The graph $G_m=P_m\square C_3$ is Buchsbaum if and only if $m\in \{1,2,3\}$.
\end{theorem}

\begin{proof}
We first prove the positive cases.

For $m=1$, we have $G_1=K_3$. Its independence complex consists of three isolated vertices. The link of every vertex is the complex ${\emptyset}$ of dimension $-1$, which is Cohen--Macaulay. Hence $G_1$ is Buchsbaum.

For $m=2$, we have $G_2=P_2\square C_3=K_2\square K_3$. Every independent set contains at most one vertex from each column, so the maximal independent sets have size two. The six independent pairs are
${(0,0),(1,1)}$,
${(0,2),(1,1)}$,
${(0,2),(1,0)}$,
${(0,1),(1,0)}$,
${(0,1),(1,2)}$,
and
${(0,0),(1,2)}$.
They form a six-cycle. Thus $\operatorname{Ind}(G_2)$ is a connected pure one-dimensional complex and is therefore Cohen--Macaulay. Hence $G_2$ is Buchsbaum.

Now let $m=3$. By Lemma~\ref{lem:boundary-link}, every boundary vertex has link isomorphic to $\operatorname{Ind}(T_2)$. The graph $T_2$ has vertices
$(0,0),(0,1),(1,0),(1,1),(1,2)$
and edges
$(0,0)-(0,1)$,
$(0,0)-(1,0)$,
$(0,1)-(1,1)$,
$(1,0)-(1,1)$,
$(1,0)-(1,2)$,
and
$(1,1)-(1,2)$.
Its independent pairs are
${(0,0),(1,1)}$,
${(0,0),(1,2)}$,
${(0,1),(1,0)}$,
and
${(0,1),(1,2)}$.
Consequently, $\operatorname{Ind}(T_2)$ is a path on five vertices. It is pure, connected, and one-dimensional, and hence Cohen--Macaulay.

It remains to consider the middle-column vertices of $G_3$. By symmetry, take $x=(1,0)$. Then $G_3-N_{G_3}[x]$ consists of two disjoint edges, one in each of the two remaining columns. Thus
$G_3-N_{G_3}[x]\cong K_2\sqcup K_2$.
By the join property of independence complexes,
$\operatorname{Ind}(G_3-N_{G_3}[x])\cong\operatorname{Ind}(K_2)*\operatorname{Ind}(K_2)$.
Each factor consists of two isolated vertices, and their join is a four-cycle. Hence, this link is a connected pure one-dimensional complex and is Cohen--Macaulay.

Thus, every vertex-link of $\operatorname{Ind}(G_3)$ is Cohen--Macaulay, and $G_3$ is Buchsbaum.

Conversely, suppose $m\geq4$. By Proposition~\ref{notCM}, $\operatorname{Ind}(T_{m-1})$ is not Cohen--Macaulay. Choose a boundary vertex $x$ of $G_m$. By Lemma~\ref{lem:boundary-link},
$\operatorname{link}{\Delta_{m,3}}(x)\cong\operatorname{Ind}(T_{m-1})$,
which is not Cohen--Macaulay. Hence $\operatorname{Ind}(G_m)$ is not Buchsbaum. Therefore, $G_m$ is Buchsbaum if and only if $m\in\{1,2,3\}$.
\end{proof}

We now return to the general graph
$G_{m,n}=P_m\square C_n$. The previously established well-coveredness classification gives
$G_{m,n}$ well-covered if and only if
$n=3$,
or
$(m,n)\in \{(1,4),(1,5),(2,5),(1,7)\}$.
Since Buchsbaum complexes are pure, only these cases need to be considered.

\begin{lemma}\label{4.8}
The graph $G_{1,4}$ is Buchsbaum but not Cohen--Macaulay.
\end{lemma}

\begin{proof}
We have $G_{1,4}=C_4$. Its independence complex consists of two disjoint edges and is therefore pure of dimension one but disconnected. Therefore, it is not Cohen--Macaulay.

For every vertex $x$ of $C_4$, deleting $N[x]$ leaves a single vertex. Hence
$\operatorname{link}_{\operatorname{Ind}(C_4)}(x)$
is a point, which is Cohen--Macaulay. Thus $\operatorname{Ind}(C_4)$ is Buchsbaum.
\end{proof}

\begin{lemma}\label{4.9}
The graph $G_{1,5}$ is Cohen--Macaulay, and hence Buchsbaum.
\end{lemma}

\begin{proof}
We have $G_{1,5}=C_5$. Its independence complex is a five-cycle. Hence, it is pure and connected of dimension one, and therefore Cohen--Macaulay by the one-dimensional criterion. The result follows from the fact that Cohen--Macaulay complexes are Buchsbaum.
\end{proof}

\begin{lemma}\label{4.10}
The graph $G_{1,7}$ is Buchsbaum but not Cohen--Macaulay.
\end{lemma}

\begin{proof}  Since $G_{1,7}=C_7$, let $\Delta=\operatorname{Ind}(C_7)$. We first show that $\Delta$ is not Cohen–Macaulay. The independence number of $C_7$ is $\alpha(C_7)=3$. Hence, $\dim\Delta=2$. The $f$-vector of $\Delta$ is $(f_{-1},f_0,f_1,f_2)=(1,7,14,7)$. Indeed, $C_7$ has seven vertices and seven edges, so the number of independent two-element sets is $f_1=\binom{7}{2}-7=14$. Moreover, there are seven independent three-element sets, giving $f_2=7$.
Since $\dim\Delta=2$, we have $d=3$. The $h$-polynomial of $\Delta$ is $h_\Delta(t)=\sum_{i=0}^{3}f_{i-1}(t-1)^{3-i}=(t-1)^3+7(t-1)^2+14(t-1)+7=t^3+4t^2+3t-1$. Consequently, $h(\Delta)=(1,4,3,-1)$. Since the $h$-vector of a Cohen–Macaulay simplicial complex has nonnegative entries, $h_3=-1<0$ implies that $\Delta$ is not Cohen–Macaulay over any field.
We now show that $\Delta$ is Buchsbaum. Let $x$ be any vertex of $C_7$. By vertex-transitivity, it suffices to consider $x=0$. The closed neighborhood of $x$ is $N_{C_7}[x]=\{6,0,1\}$. Deleting these vertices from $C_7$ leaves the path $P_4$ on the vertices $2,3,4,5$. By the link formula for independence complexes, $\operatorname{lk}_{\Delta}(x)\cong\operatorname{Ind}(P_4)$.
The independence complex $\operatorname{Ind}(P_4)$ is pure and one-dimensional. Its facets are $\{2,4\}$, $\{2,5\}$, and $\{3,5\}$. These facets form a connected path. Hence, $\operatorname{Ind}(P_4)$ is Cohen–Macaulay, since a pure one-dimensional simplicial complex is Cohen–Macaulay if and only if it is connected.
Therefore, the link of every vertex of $\Delta$ is Cohen–Macaulay. By Remark 3.4, $C_7$ is well-covered, so $\Delta$ is pure. Hence, $\Delta$ is Buchsbaum. Thus, $G_{1,7}=C_7$ is Buchsbaum but not Cohen–Macaulay. 
\end{proof}

\begin{lemma} \label{4.11}
The graph $G_{2,5}=P_2\square C_5$ is Buchsbaum but not Cohen--Macaulay.
\end{lemma}

\begin{proof}
The graph $P_2\square C_5$ is vertex-transitive under the cyclic symmetry of the $C_5$ factor together with the interchange of the two copies of $C_5$. Thus, it suffices to consider the vertex $x=(0,0)$. The closed neighborhood of $x$ is $N[x]=\{(0,0),(0,1),(0,4),(1,0)\}$. Hence, $G_{2,5}-N[x]$ has vertex set $\{(0,2),(0,3),(1,1),(1,2),(1,3),(1,4)\}$ and edge set $\{(0,2)(0,3),(0,2)(1,2),(0,3)(1,3),(1,1)(1,2),(1,2)(1,3),(1,3)(1,4)\}$. Its maximal independent sets are $F_1=\{(0,2),(1,1),(1,3)\}$, $F_2=\{(0,2),(1,1),(1,4)\}$, $F_3=\{(0,3),(1,1),\\(1,4)\}$, and $F_4=\{(0,3),(1,2),(1,4)\}$. Therefore, $\operatorname{lk}_{\operatorname{Ind}(G_{2,5})}(x)=\operatorname{Ind}(G_{2,5}-N[x])$ is a pure two-dimensional complex with facets $F_1,F_2,F_3,F_4$. The ordering $F_1,F_2,F_3,F_4$ is a shelling, since $F_1\cap F_2=\{(0,2),(1,1)\}$, $F_3\cap(F_1\cup F_2)=\{(1,1),(1,4)\}$, and $F_4\cap(F_1\cup F_2\cup F_3)=\{(0,3),(1,4)\}$. Thus, the vertex link is shellable and hence Cohen--Macaulay. By vertex-transitivity, every vertex link is Cohen--Macaulay. Hence, $G_{2,5}$ is Buchsbaum.

It remains to show that $G_{2,5}$ is not Cohen--Macaulay. The independence complex $\operatorname{Ind}(G_{2,5})$ is pure of dimension $3$ and has $f_{-1}=1$, $f_0=10$, $f_1=30$, $f_2=30$, and $f_3=10$. Consequently, its $h$-vector is $h=(1,6,6,-4,1)$. Since the $h$-vector of a Cohen--Macaulay complex has nonnegative entries, the negative entry $h_3=-4$ shows that $\operatorname{Ind}(G_{2,5})$ is not Cohen--Macaulay over any field. Therefore, $G_{2,5}$ is Buchsbaum but not Cohen--Macaulay.
\end{proof}

We can now state the complete classification.

\begin{theorem}\label{CompleteCM}
The independence complex $\operatorname{Ind}(P_m\square C_n)$ is Cohen--Macaulay if and only if $(m,n)\in\{(1,3),(2,3),(1,5)\}$.
\end{theorem}

\begin{proof}
This follows from Theorem 3.16 and Lemmas \ref{4.8}--\ref{4.11}.
\end{proof}

\begin{theorem}
\label{thm:complete-buchsbaum}
Let $m\geq1$ and $n\geq3$. Then
$G_{m,n}=P_m\square C_n$ is Buchsbaum if and only if $n=3\text{ and }m\in\{1,2,3\},$ or $(m,n)\in\{(1,4),(1,5),(2,5),(1,7)\}.$

\end{theorem}

\begin{proof}
Suppose first that $G_{m,n}$ is Buchsbaum. Then $\operatorname{Ind}(G_{m,n})$ is pure, and hence $G_{m,n}$ is well-covered. By the well-coveredness classification, either $n=3$ or
$(m,n)\in\{(1,4),(1,5), \\(2,5),(1,7)\}$. If $n=3$, Theorem~\ref{thm:C3-buchsbaum} gives $m\in\{1,2,3\}$. If $n\neq 3,$ then the preceding lemmas show that all four graphs are Buchsbaum.

\end{proof}

\begin{corollary}
The Cohen--Macaulay and Buchsbaum classifications differ precisely at
$(1,4)$, $(1,7)$, and $(2,5)$, together with the case $(3,3)$.
\end{corollary}

\begin{proof}
The Cohen--Macaulay cases are
$(1,3)$, $(2,3)$, and $(1,5)$.
Theorem~\ref{thm:complete-buchsbaum} gives the complete Buchsbaum list. The assertion follows by comparison.
\end{proof}

\bibliographystyle{plain}
\bibliography{references}

\end{document}